
\magnification=\magstephalf
\baselineskip=13pt
\hsize = 5.5truein
\hoffset = 0.5truein
\vsize = 8.25truein
\voffset = 0.2truein

\tolerance=10000
\overfullrule=0pt



\input epsf

\def\draw#1#2#3{
  \removelastskip
  \vskip 0cm plus #1cm
  \goodbreak
  \vskip 0cm plus -#1cm
  \vskip 0.3cm
  \dimen0=#2 true bp
  \ifdim\dimen0>\hsize\dimen0=\hsize\fi
  \epsfxsize=\dimen0
  \centerline{\epsfbox{#3.eps}}
  \count255=\epsfurx
  \advance\count255 by-\epsfllx
  \message{[#3.eps bp width = \number\count255]}
  \goodbreak
  \vskip 0.3cm
}





\def\bpad{\kern 0.15em}

\def\bmatrix#1{\left[\,\vcenter{\bpad\halign{
  &\bpad\hfil$##$\bpad\hfil\cr
  #1}\bpad}\,\right]}

\def\rightarrowmat#1{
  \setbox1=\hbox{\bpad$\scriptstyle\bmatrix{#1}$\bpad\bpad}
  \mathop{\hbox to\wd1{\rightarrowfill}}\limits^{\box1}}

\def\leftarrowmat#1{
  \setbox1=\hbox{\bpad\bpad$\scriptstyle\bmatrix{#1}$\bpad}
  \mathop{\hbox to\wd1{\leftarrowfill}}\limits^{\box1}}


\def\Box{
\hskip 12pt minus 4pt
\vrule height 1.2ex width .9ex depth .1ex
\par \goodbreak \vskip .3 cm
}

\centerline{\bf  MONOMIAL RESOLUTIONS }

\vglue .3cm

\centerline{\bf Dave Bayer \ \  Irena Peeva  \ \   Bernd Sturmfels }

\vglue .3cm

\beginsection  1. Introduction

Let $M$ be a monomial ideal in the polynomial ring $S=k[x_1,\dots ,x_n]$
over a field $k$. We are interested in the problem of resolving $S/M$ over $S$.
The difficulty in resolving minimally is reflected in the fact that the
homology of arbitrary simplicial complexes can be encoded (via the 
Stanley-Reisner correspondence) into the multigraded Betti numbers of $S/M$, 
[St]. In particular, the minimal free resolution may
depend on the characteristic of $k$. In the study of
monomial resolutions  the following main directions have been pursued:
construction of the minimal free resolution for special types of ideals
[EK], construction of  structured non-minimal resolutions [Ly],
bounds for the numerical invariants
of the minimal free resolution [Bi,Hu], algorithms 
and software for computing resolutions and  Hilbert functions [BS].

We present an  approach  for resolving  $S/M$ by
 encoding the entire resolution into a single simplicial complex;
the obtained resolution is minimal generically.

\vskip .2cm

\noindent {\bf Construction~1.1.}
{\sl (Monomial resolution from a labeled simplicial complex) }
\hfill \break
Let $\Delta $ be a simplicial complex whose vertices are labeled by
generators of $M$. We label each face of $\Delta$ by the least common
multiple of its vertices. The exponent vectors of these monomials
define an ${\bf N}^n$-grading of $\Delta$.
Let ${\bf F}_{\Delta}$  be the ${\bf N}^n$-graded  chain complex of
$\Delta$ over $S$. It is obtained from the usual
chain complex by homogenizing the differential. For example,
an edge  $\{1,2\}$, whose vertices are labeled by $x^3y^4$ and $xz^3$, is
mapped by the differential to $z^3\{1\}-x^2y^4\{2\}$.
If the complex   ${\bf F}_{\Delta}$ is exact then we call it the
{\it resolution defined by the labeled simplicial complex $\Delta$}.
Such a resolution
is characteristic-free and a DG-algebra
(associative commutative differential graded algebra).
In this case the ${\bf N}^n$-graded Hilbert series of $S/M$ equals
the ${\bf N}^n$-graded  Euler characteristic of $\Delta $
divided by $(1-x_1)\cdots (1-x_n)$.
\Box

\noindent{\bf Example~1.2.}
The classical example for Construction 1.1
is {\sl Taylor's resolution}, which
arises when  $\Delta$ is the full simplex
on the minimal generators of $M$.
An algebraic review of Taylor's resolution
and a short proof of its exactness are given in Section 2.
For a large number of generators this resolution
is very far from minimal and it is inefficient  for obtaining
bounds on Betti numbers and regularity.
\Box

\noindent{\bf Definition 1.3.}
A monomial ideal $M$ is called
{\it generic} if no variable $x_i$ appears with the same
non-zero exponent in two distinct minimal generators of $M$.

\vskip .2 cm

Almost all monomial ideals are generic, in the sense that those
which fail to be generic  lie on finitely many hyperplanes
in the matrix space of exponents.
In this paper we prove that the minimal
free resolution for any generic monomial ideal $M$
comes from a simplicial complex $\Delta_M$, which we call the
{\it Scarf complex} of $M$.

\vskip .2 cm
\noindent{\bf Example~1.4.}
The Scarf complex of $\, M\, =\, \langle \,
x^2z^3,\, x^3z^2,\, xyz,\, y^2\, \rangle \,$  is
a triangle connected to an edge:

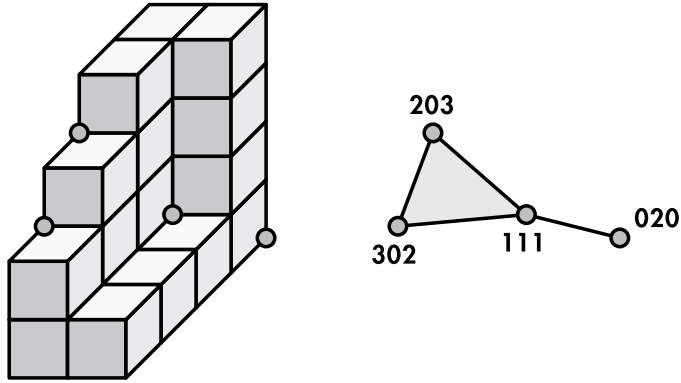
\draw4{254}{figA} 

\noindent
The triangle is labeled by $x^3yz^3$, the edges of the triangle are
labeled by $\,x^3z^3, \, x^2yz^3, \,x^3yz^2$,
and the other edge by $xy^2z$.
The minimal free resolution  of $S/M$ is the
${\bf N}^3$-graded chain complex of
the simplicial complex depicted above:
 $$
 0 \rightarrow S
 \rightarrowmat{y \cr -x \cr z \cr 0 \cr}
 S^4
 \rightarrowmat{-x & -y & 0 & 0 \cr z & 0 & -y & 0 \cr
 0 & xz^2 & x^2z & -y \cr 0 & 0 & 0 & xz \cr}
 S^4
 \rightarrowmat{x^2z^3 &\, x^3z^2 & \, xyz & \, y^2 \cr}
 S \rightarrow S/M \rightarrow 0.
 $$
The Hilbert series of $S/M$ equals
$$\sum \{ \, x^a y^b z^c \,: \,x^a y^b z^c \not\in M \,\}
\quad =  \quad {P(x,y,z)\over (1-x)(1-y)(1-z)}, \qquad \hbox{ where}$$
$$ P(x,y,z) \; = \;
 1-x^2z^3-x^3z^2-xyz-y^2+x^3z^3+x^2yz^3+x^3yz^2+xy^2z-x^3yz^3
 .$$
This polynomial is the ${\bf N}^3$-graded
Euler characteristic of the Scarf complex.
\Box

 An algebraic construction of the Scarf complex is
given in Section~3.
In Section~4 we construct a convex polytope whose
boundary naturally contains the Scarf complex.
This is based on work of Herbert Scarf  in
mathematical economics [Sca].
In the artinian case the Scarf complex is a regular
triangulation of a simplex. In general it
need not be pure (Example 1.4), and it need not be shellable (Example 4.2),
but it is always contractible (Theorem 4.1).
In Section 7 we raise the question of characterizing all Scarf 
complexes of monomial ideals in $n$ variables.
This is related to the {\it order dimension} of posets.
In  Section~8 we relate
the Scarf complex to the irreducible decomposition of $M$
and to the Cohen-Macaulay property.

If $M$ is a monomial ideal which is not generic, then typically the
minimal free resolution of $S/M$ does not come from Construction 1.1.
Two notable obstructions are dependence on the characteristic of $k$
and non-existence of a DG-algebra structure.
In Section 5 we obtain a {\sl non-minimal} free resolution
which comes from a simplicial complex and
has length at most the number of variables.
Our construction is based on  {\it degeneration of exponents}:
we deform $M$ to a nearby generic monomial ideal
$M'$ by using monomials with real exponents and moving the exponents of the
generators
to generic position, we compute the
Scarf complex $\Delta_{M'}$ of the generic ideal $M'$, and then
we label the vertices of $\Delta_{M'}$
with the generators of the original ideal $M$. This defines
a non-minimal resolution for $S/M$. As an application, in
Section~6 we show that the Betti numbers satisfy
the inequalities of the Upper Bound Theorem for Convex Polytopes.

\beginsection  2.   Taylor's resolution

Let $M=\langle m_1, \dots , m_r\rangle$ be a
monomial ideal in $S=k[x_1,\dots ,x_n]$.
For a subset $I$ of $\{1,\dots ,r\}$ we set
$m_I:=lcm(m_i \  \vert \  i\in I)$. Let
$a_I \in {\bf N}^n$ be the exponent vector of $m_I$ and
 $S(-a_I)$ the free $S$-module with one generator in multidegree $a_I$. The
{\it   Taylor resolution} of $S/M$ is the ${\bf Z}^n$-graded module ${\bf F}=
{\bigoplus}_{I\subseteq \{1,\dots ,r\} } S(-a_I)$
 with basis denoted by
${\{e_I\}}_{I\subseteq \{1,\dots ,r\} }$
and equipped with the differential
$$d(e_I)\ = \ {\sum }_{i\in I} sign(i,I)\cdot
{m_I \over m_{I \setminus  i}} \cdot e_{I \setminus  i}\ ,\eqno(2.1)$$
where $sign(i,I)$ is $(-1)^{j+1}$ if $i$ is the $j$th element in the
ordering of $ I$. This is a free resolution of $S/M$ over $S$
having length $r$ and $2^r$ terms.
It is very far from minimal if $r\gg n$. 
A smaller resolution based on Taylor's
resolution was constructed in  [Ly].

Every simplicial complex $\Delta$ on $\{1,\dots ,r\}$
defines a submodule $\,{\bf F}_{{\Delta}} :=
{\bigoplus}_{I\in {\Delta} } S(-a_I)$ of the Taylor resolution
${\bf F}$ which is closed under the differential (2.1).
The complex ${\bf F}_\Delta$ is the one
introduced in Construction 1.1, with
each face $ I $ labeled by $m_I$.
We now use reduced simplicial
homology to determine when ${\bf F}_\Delta$ is exact.

\proclaim Lemma~2.1. The complex ${\bf F}_{{\Delta}}$ is exact
if and only if for every monomial $m$ the simplicial complex
${\Delta}[m]=\{I\in \Delta  \,\vert\,  m_I \hbox{\ divides\ }m\}$
 is empty or acyclic over $k$.

\noindent {\sl Proof:}
Since ${\bf F}_{{\Delta}}$ is  ${\bf N}^n$-graded
it suffices to check exactness in each multidegree.
The component of ${\bf F}_{{\Delta}}$ in multidegree $m$
is a complex of finite-dimensional $k$-vector spaces, which
can be identified with the
chain complex of ${\Delta}[m]$ over $k$.
\Box

\noindent {\bf Remark~2.2.}
Lemma~2.1 shows that Taylor's complex ${\bf F}$ itself is exact.
This is the case when ${\Delta}$ is the full $(r-1)$-simplex; cf. Example~1.2.
For each multidegree $m$  we have that ${\Delta}[m]$
is the full simplex on the minimal generators which divide $m$.
\Box

\beginsection 3. Generic monomial ideals

For any monomial ideal $ M = \langle m_1, \ldots, m_r \rangle $ 
we define a simplicial complex:
 $${\Delta}_M \; := \;
 \bigl\{\,I \subseteq \{1,\ldots,r\} \,\,
 \vert \,\,\ m_I\not= m_J \ \,\, \hbox{\rm for all  }\ J \subseteq
\{1,\ldots,r\} \, \hbox{ other than }  \, I \,\bigr\}  \eqno (3.1) $$
This complex was first introduced in
1973 by Herbert~Scarf [Sca] in the context of mathematical economics.
For this reason we call
$\Delta_M$ the {\it Scarf complex} of $M$.

 Taylor's resolution ${\bf F}$
is a direct sum of the minimal free resolution of $S/M$
and trivial complexes
$ 0 \longrightarrow S(-a_I) \longrightarrow S(-a_I) \longrightarrow 0$. On
the other hand, if $I\in {\Delta}_M$ then  ${\bf F}$ has a unique minimal
generator in multidegree $a_I$. Therefore, the minimal free resolution of
$S/M$ always contains the complex  ${\bf F}_{{\Delta}_M}$, but is larger in
general.  For every monomial ideal $M$ and every face $I$ of $\Delta_M$,
the minimal free resolution of $S/M$ has a unique generator
in multidegree $a_I$.

\proclaim Lemma~3.1. If all non-zero Betti numbers of $S/M$ are
concentrated in the multidegrees $a_I$
 of the faces $I$ of ${\Delta}_M$,
then ${\bf F}_{{\Delta}_M}$ is the minimal free resolution of $S/M$.

\noindent {\sl Proof: }
 If the minimal free
resolution  of $S/M$ is strictly larger than  ${\bf F}_{{\Delta}_M}$,
then the Taylor resolution has at least
 two basis elements in some multidegree
$a_I$ for $I \in {\Delta}_M$.
This contradicts the definition of $\Delta_M$.
\Box

We now consider generic monomial ideals (see Definition 1.3).

\proclaim Theorem~3.2.
Let $M$ be a generic monomial ideal.
Then the complex  ${\bf F}_{{\Delta}_M}$
defined by the Scarf complex $\Delta_M$
is a minimal free resolution of $S/M$ over $S$.

\noindent {\sl Proof: }
If $ I \in \Delta_M$ and $i \in I$ then
$m_{I \setminus i}$ properly divides $m_I$. 
Thus we see directly from (2.1) that ${\bf F}_{{\Delta}_M}$
is minimal. It remains to show that ${\bf F}_{{\Delta}_M}$ is exact.
Consider any multidegree $a_I$ with $I \not\in \Delta_M $.
Choosing $I$ minimal with respect to inclusion, we
may assume  $\, m_I = m_{I \cup i} \,$ for some
$i \in \{1,\ldots,r\}\backslash I $.
The monomials $m_i$ and $m_I$ have different exponents in any fixed  variable
because $M$ is generic. By Hochster's formula [Ho, Theorem 5.1], 
the $j$th Betti number in multidegree $a_I$
equals the dimension of the reduced $(j-2)$nd
homology of the simplicial complex on $\{x_1,\ldots,x_n\}$
generated by the supports of $\,m_I / m_s $, where $m_s$
runs over all generators which divide $m_I$.
Taking $s=i$, we see that this simplicial complex
is the full simplex on the support of $m_I$,
hence contractible. Hence the
$j$th Betti number in multidegree $a_I$
is zero unless $I \in \Delta_M$.
Now apply Lemma 3.1.
\Box

\proclaim Corollary~3.3. Let $M$ be a generic monomial ideal.
\item{(1)}
The number of $j$-faces of the Scarf complex $\Delta_M$
equals the total Betti number
$\,\beta_{j+1}(S/M)\, =\, dim_k \,Tor_{j+1}^S ( S/M, k)$.
\item{(2)} The minimal free resolution of $S/M$
is characteristic free.  It is ${\bf N}^n$-graded and in
each multidegree the Betti number is either $0$ or $1$.
\item{(3)}
The ${\bf N}^n$-graded Hilbert series of $S/M$
(i.e.~the sum over all monomials not in $M$) is
$$
{ {\sum}_{ I\in {\Delta}_M }(-1)^{|I|} \cdot {m_I} \over
(1- x_1)\cdots (1-x_n)}, $$
and there are no cancellations in the alternating sum in the numerator.

The multigraded Betti numbers of  any monomial ideal
determine its Hilbert series. The converse may fail.
Theorem 3.2 and Corollary 3.3 show that, for a
generic monomial ideal,  the multigraded Hilbert series
determines the minimal free resolution.
In particular, they are equivalent from 
a computational point of view. 

\vskip .2cm

\proclaim Corollary~3.4. If $M$ is a generic monomial ideal
then the minimal free resolution of $S/M$ is a DG-algebra.

We first recall the necessary definition:
Let $\bf A$ be a resolution of a cyclic
$S$-module $S/L $ with differential $d$ and suppose
 $\bf A$ is a graded algebra with multiplication $*$.
Then  $({\bf A},d,*)$ is called a {\it commutative associative
differential graded $S$-algebra (=DG-algebra)}
if $*$ is associative and for any homogeneous elements
$\alpha,\beta \in \bf A$ the product is skew commutative,
$\alpha * \beta =(-1)^{deg(\alpha ) deg(\beta )}\beta * \alpha$,
and the  Leibniz rule $d(\alpha *\beta)=d(\alpha)*\beta +(-1)^{deg(\alpha
)}\alpha *d(\beta)$ holds.
In this case we say that $\bf A$ has a DG-algebra structure.
The classical example
is the Koszul complex resolving a regular sequence.
Advantages of the existence of such a structure are that it
yields lower bounds on the Betti numbers of $S/L$
and  information on the infinite minimal free resolution of $k$ over $S/L$. 
 The   Taylor resolution of $S/M$ for $M$ a monomial ideal has
a DG-algebra structure, however the minimal free resolution might not:

\vskip .2cm

\noindent {\bf Example~3.5.} {\sl (Backelin)}
For $M =\langle x^2,xy^2z,y^2z^2,yz^2w,w^2\rangle$
the minimal free resolution of $S/M$ has no DG-algebra structure.
This was shown in [Av, 5.2.3]. \Box

\vskip .2 cm

Every resolution admits a commutative differential graded structure.
The problem is to ensure associativity. This usually requires a complex
way of multiplying; cf.~the DG-algebra structure 
in [Sr] on the minimal resolution of  $S/\langle x_1,\dots ,x_n\rangle^p $.
A notable advantage of Construction~1.1 is that
leads to a simple multiplication:

\vskip .1cm

\noindent {\sl Proof of Corollary 3.4: } We define a multiplication on
 ${\bf F}_{{\Delta}_M} $ by :
$$
e_I*e_J\ :=\  \left\{
\eqalign {
{\ sign (I,J) \cdot {m_Im_J \over m_{I\cup J}}\cdot
 e_{I\cup J}} \ &\quad{\rm if}\ I\cup J \in {\Delta}_M
\,\,\, {\rm and} \, \,\, I\cap J=\emptyset ,\ \cr
\ 0\quad \quad &\quad{\rm otherwise.}
}
\right.
$$
Here $sign (I,J)$ is  defined as follows: let $I=\{i_1<\dots <i_q\}$ and
$J=\{j_1<\dots <j_p\}$, and $\tau$ be the permutation
which permutes $i_1,\dots ,i_q,
j_1,\dots ,j_p$ into a sequence of increasing numbers, then $sign(I,J)$ 
is the sign of $\tau$. Note that this multiplication is multigraded. 
Associativity and skew commutativity are obvious.
Straightforward computations verify that the
Leibniz rule holds. \Box

\beginsection  4. Convexity

Polytopes are a powerful tool for structuring
combinatorial data appearing in algebra and algebraic 
geometry. For example, Newton polytopes play a 
significant role in computer algebra, 
singularity theory and toric geometry.
In this section we present the polytope
underlying the minimal free resolution of a generic monomial ideal $M$.

Let $M[m]$ be the subideal of $M$
generated by the minimal generators of $M$ which
divide a given monomial $m$.  Then $M[m]$ is a generic
monomial ideal as well, and its Scarf complex equals
$$ {\Delta}_{M[m]} \quad = \quad {\Delta}_M[m]. \eqno (4.1) $$
Lemma~2.1 and  (4.1) imply that the exactness of
${\bf F}_{{\Delta}_M}$ for all $M$
is equivalent to the acyclicity of ${\Delta}_M$
for all $M$. In fact, the following stronger result holds.

\proclaim Theorem~4.1. The  Scarf complex ${\Delta}_M$
of a generic monomial ideal $M$  is contractible.

The Scarf complex need not be pure, and
it need not be shellable either:

\vskip .2cm

\noindent
{\bf Example~4.2.}
Consider the generic monomial ideal
$$ M \quad := \quad
\langle \,
xyz, \,
x^4y^3, \,
x^3y^5, \,
y^4z^3, \,
y^2z^4, \,
x^2z^2 \,
\rangle .$$
The Scarf complex consists of two triangles and an edge meeting
at a vertex.

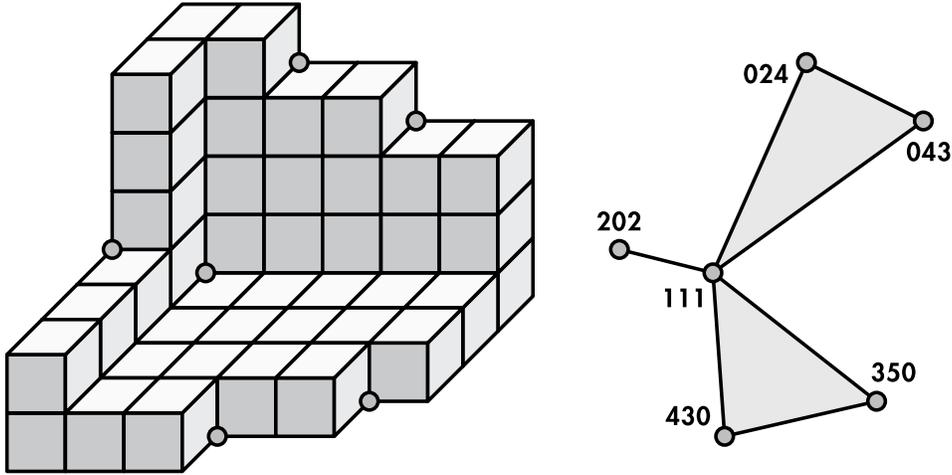
\draw4{357}{figB} 

\noindent
Here $\Delta_M$ is contractible, but not shellable,
even in the non-pure sense of [BW].
\Box

Theorem~4.1 will be derived from the following convexity theorem,
which is an extension  of Theorem 2.8.4 in [Sca, \S 2.8].
Scarf calls the faces of $\Delta_M$ ``primitive sets''
and relates them to dual feasible bases in linear programming.

\proclaim Theorem~4.3. {\rm (Scarf 1973) }
Let $M$ be a generic monomial ideal. There exists a
polytope $P_M$ with $r$ vertices in ${\bf R}^n$ such that
${\Delta}_M$ is isomorphic to
the subcomplex of the boundary of $P_M$ consisting of all
faces supported by a strictly positive inner normal vector.

\noindent {\bf Construction~4.4. } One possible choice
of a polytope $P_M$ satisfying Theorem~4.3 is as follows.
Let $a_i = (a_{i1}, a_{i2}, \ldots,a_{in}) \in {\bf N}^n$ be
the exponent vector of the $i$-th minimal generator of $M$.
Scarf's original definition of the Scarf complex is as follows:
$$  I \in \Delta_M \quad \iff \quad
 \forall \,i \in \{1,\ldots,r\} \,\, \exists \, j \in \{1,\ldots,n \} \,:
\,\,\,   a_{I j} \leq a_{ij} \eqno (4.2) $$
where $a_{Ij} = max \{\, a_{ij} \,| \, i \in I  \,\}$.
We fix a sufficiently large real number $t \gg 0$ and define
$P_M$ as the convex hull of the point set
$$ \{\, ( a_{i1}^t, a_{i2}^t,\dots ,
a_{in}^t )\  \, \vert \,\, \ 1\le i\le r\ \}
\quad \subset \quad {\bf R}^n . \eqno (4.3) $$
The combinatorial type of $P_M$ is independent of $t$ for  large $t$.

We remark that there are other ways of constructing $P_M$;
for example, take the convex hull of the points
$\{\ ( t^{a_{i1}},\dots ,t^{a_{in}})\ \vert \ 1\le i\le r\ \}\,$
for sufficiently large $t$.
\Box

\noindent {\sl Proof of Theorem~4.3: }
We identify each face of $P_M$ with a subset
$I \subseteq \{1,\ldots,r \}$, namely, the
indices of vertices  which lie on that face.
Let $I$ be a face of $P_M$ with
inner normal vector $(w_1,\ldots,w_n) $ where
$w_j > 0$ for $ j=1,\ldots,n$.
We may assume $$ \forall \, i \in
\{1,\ldots,r \} \,: \,\,\,\, w_1 \cdot a_{i1}^t \, +\,
\cdots \, +\,  w_n \cdot
a_{in}^t \,\, \geq \,\, 1  , \eqno (4.4) $$
and equality holds in (4.4) if and only if $i \in I $.
This implies
$$ \forall \, j \in \{1,\ldots,n\} \,\, : \quad
a_{Ij}
\,\,\, \leq \,\,\, w_j^{- 1/t} . \eqno (4.5) $$
At least one of the summands in (4.4)  is greater or equal to $1/n$.
This implies
$$
\forall \,i \in \{1,\ldots,r \}
\,\,\,\exists \,j \in \{1,\ldots,n \} \,: \quad
a_{ij} \,\geq\,
w_j^{-1/t} \cdot n^{- 1/t} \,
\geq \,
a_{Ij} \cdot n^{- 1/t} . $$
Now let $t \rightarrow \infty$.
Then we can erase the factor
$\,n^{- 1/t} \rightarrow 1 $, and (4.2) is satisfied.

\vskip .2cm

For the converse we consider the special case where $M$ is artinian.
We first show that $\Delta_M$ is pure of dimension $n-1$.
We may assume that
$ \,a_{ij} = 0 \,$ for
$1 \leq i \not= j \leq n$ and
$\, a_{li} < a_{ii} \,$
for $1 \leq i \leq n < l \leq r $.
Consider any $I \in \Delta_M$ with $| I | < n$.
Then there exists $i \in I$ such that
$a_i$ and $a_I$ agree in at least two coordinates,
say $j$ and $j'$.
Consider the set
$$ {\cal S} \quad := \quad \bigl\{
\, l \in \{1,\ldots,r\} \,\,\,| \,\,\,
\,a_{lj} > a_{ij} = a_{Ij} \,\, \hbox{and} \,\,
\,\forall i \in I \,\, \exists s \,: a_{is} > a_{ls}\, \bigr\}. $$
It is non-empty since $j \in {\cal S}$.
Select $l \in {\cal S}$ with $a_{lj}$ smallest.
Then $\,I \, \cup \, \{l\} \,\in \,\Delta_M$.

We next consider the
{\it oriented matroid} [BLSWZ] of the configuration
in (4.3) plus the origin. Set
$\,a_{0j} := 0\,$ for all $j\,$ and define for $0 \leq i_0 <  \ldots < i_n
\leq  r$
 $$ [ i_0,i_1,\ldots,i_n ] \,\,\, := \,\,\,
  sign \bigl( det  \pmatrix{
  1 &  a_{i_0 1}^t &  a_{i_0 2}^t & \cdots &  a_{i_0 n}^t \cr
  1 &  a_{i_1 1}^t &  a_{i_1 2}^t & \cdots &  a_{i_1 n}^t \cr
  \vdots &   \vdots &   \vdots & \ddots & \vdots \cr
  1 &  a_{i_n 1}^t &  a_{i_n 2}^t & \cdots &  a_{i_n n}^t \cr} \bigr). $$
Let $I = \{i_1,\dots,i_n\}$ be a maximal face of $\Delta_M$.
There exists a unique permutation
$\sigma = (\sigma_1,\ldots,\sigma_n)$  of $\,I \,$
such that $\, a_{Ij} = a_{\sigma_j j} \, $ for all $j$.
Since $t \gg 0$, we have
$$ sign(\sigma)  \,\,= \,\,    [\, 0, i_1,\ldots,i_n ]. $$
For $j \not\in I$ there exists an $s$ with
$a_{js} > a_{\sigma_j s}$. Laplace expansion
along the $s$-th column gives
$$ [\, j, i_1,\ldots,i_n ] \,\,= \,\, -sign(\sigma)  \quad
\hbox{for all} \,\, j \not\in I  . $$
This shows that $I$ is a facet of $P_M$ which is visible from the
origin. Since $P_M$ intersects each coordinate axis, 
the normal vector of $I$ is strictly positive.

Now drop the assumption that $M$ is artinian
and let $I$ be a maximal face of $\Delta_M$.
Let $M'$ be an artinian ideal obtained from $M$
by adding large powers of the variables.
Then $I$ lies in $\Delta_{M'}$ as well.
Therefore $I$ is a face of $P_{M'}$  having a
positive inner normal vector, and since $P_{M} \subseteq P_{M'}$,
that same positive vector is minimized over $P_M$ at $I$.
\Box

It follows from our discussion that
the Scarf complex $\Delta_M$ is pure and shellable when
$M$ is artinian, that is, when every
variable $x_i$ appears to some power in $M$.

\proclaim Corollary~4.5. If $M$ is artinian,
then ${\Delta}_M$ is a regular triangulation of the $(n-1)$-simplex.

\noindent {\sl Proof: }
 $P_M$ lies in the positive orthant
and intersects each coordinate axis.
Each face of $P_M$ visible from the origin
is a simplex. The set of these faces is $\Delta_M$.
\Box

The definition of {\it regular triangulations}
and their basic properties can be found in [Lee].
We shall now prove the main result of this section:

\vskip .2cm

\noindent {\sl Proof of Theorem~4.1: }
Let $I$ be a face of $P_M$ and
let $\,{\cal N}(I) \,$ denote the inner normal
cone of $P_M$ at $I$.
By Theorem 4.3, the Scarf complex $\Delta_M$ consists of all
faces $I$ such that ${\cal N}(I)$ intersects the open positive
orthant. Choose an $\epsilon > 0 $ such that
$I \in \Delta_M $ if and only if
the closed cone $\,\overline{{\cal N}(I)}\, = \,
 \bigcap_{i \in I }  \overline{{\cal N}(\{i\})}\,$ intersects the
$(n-1)$-simplex
$$ T \quad :=\quad
\bigl\{ \, (u_1,\ldots,u_n) \in {\bf R}^n \,:\,
u_1 + \cdots + u_n =1 \,, \, u_i \geq \epsilon \, \,\,
\hbox{for all $i$} \,\,\bigr\}. $$
Then $\,\, \bigcap_{i \in I } \bigl(\, \overline{{\cal N}(i)} \,\cap \,
T \, \bigr)  \,$ is non-empty if and only if $I \in \Delta_M$.
Thus $\, \bigl\{
\overline{{\cal N}(i)} \,\cap \, T \,  \bigr\}_{1 \leq i \leq r} \,$
is a cover of $T$ by polytopes. The nerve
of this cover equals $\Delta_M$. Using
Borsuk's Nerve Lemma (cf.~[Bj, Thm.~10.6]),
we see that $\Delta_M$ is homotopy equivalent to $T$.
\hfil\break\hbox{}\hfill\Box

The polytope $P_M$ provides a geometric construction
of the Scarf complex $\Delta_M$ and hence
of the minimal free resolution ${\bf F}_{\Delta_M}$  of $S/M$.

\vskip .2cm

\noindent
{\bf Example~4.6.}
{\sl Every bivariate monomial ideal is generic} and can be written as
$$ M \quad = \quad
\langle \,
x^{a_1} y^{b_1} ,\,
x^{a_2} y^{b_2} ,\,\ldots, \,
x^{a_r} y^{b_r} \,\rangle,$$
where $a_1 > a_2 > \cdots > a_r$ and $b_1 < b_2 < \cdots < b_r$.
For $t \gg 0 $  the polytope
$\,P_M \, = \, conv \,\{ \,(a_i^t,b_i^t ) \,\,|\,\, i=1,\ldots, r \,\}\,$
is a convex $r$-gon. The Scarf complex $\Delta_M$ consists of the
$r-1$ edges of $P_M$ which are visible from the origin.
This triangulation of the $1$-simplex  is a
visualization of the minimal free resolution of $S/M$:
$$  {\bf F}_{\Delta_M} \; : \quad
 0 \; \longrightarrow \;
\bigoplus_{i=1}^{r-1}  \, S \bigl( -(a_i,b_{i+1})  \,\bigr)
 \; \longrightarrow \;
\bigoplus_{j=1}^{r}  \,S \bigr(-(a_j,b_j) \,\bigr)
 \; \longrightarrow \; S .\  $$

\noindent
{\bf Example~4.7. } Consider the ideal
$ M\, =\,\langle \, x^4, \, y^4, \,\allowbreak z^4, \,\allowbreak x y^2
z^3, \,\allowbreak x^3 y z^2, \,\allowbreak x^2 y^3 z \rangle $.
 This ideal is artinian and generic.
 Label the generators $1,2,3,4,5,6$ in the order they are listed. The
Scarf complex of $M$ consists of the $7$ triangles $$
 \{1,2,6\},\{1,3,5\},\{1,5,6\},\{2,3,4\},\{2,4,6\},\{3,4,5\},\{4,5,6\}.$$
This is the boundary complex of an octahedron  with the one facet
$\{1,2,3\}$ removed. The ``staircase diagram'' of $M$ 
and the triangulation $\Delta_M$ are shown below:

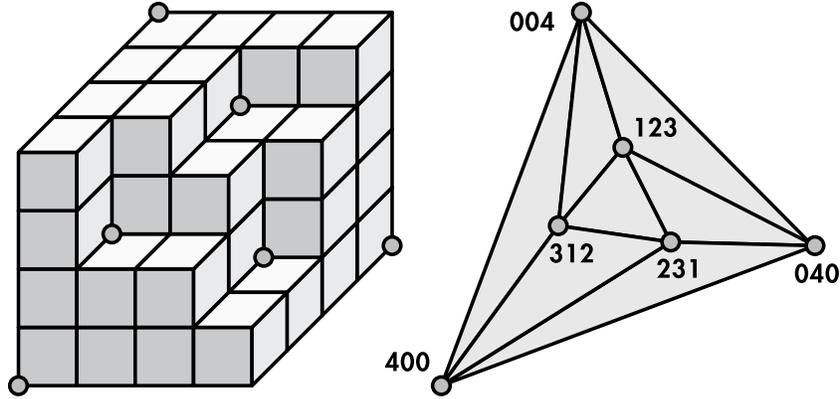
\draw4{314}{figC} 

\noindent
The minimal free resolution of $S/M$ is derived from  $\Delta_M$
by Construction~1.1. \Box

\noindent
{\bf Remark~4.8. }
In Section~3 we proved that the multigraded Betti numbers
for generic $M$ are $1$ at each corner and $0$ elsewhere.
Thus, the faces of $\Delta_M$ are labeled by the corners
in the staircase diagram of $M$; see the pictures above. 
\Box

\beginsection  5. Degeneration

Let $M = \langle m_1,\ldots,m_r \rangle $ be an arbitrary monomial ideal;
say, $M$ is not generic. In this section we construct a
(generally non-minimal) free resolution of $S/M$
by degenerating the exponent vectors of the generators of $M$.
This approach has the following advantages:
\item{$\bullet$} The resolution by degeneration
has length at most the number of variables; thus
in general it is much smaller and shorter than   Taylor's resolution.
\item{$\bullet$} This resolution is a DG-algebra \
(by the same argument as in Corollary 3.4).
\item{$\bullet$} Section 4
relates the Betti numbers to the $f$-vector of a polytope.

\vskip .2cm

\noindent {\bf Construction~5.1.}
Let $\,\{\ a_i=(a_{i1},\dots ,a_{in})\ \vert \ 1\le i\le r\ \}$ be the
exponents of the minimal generators of $M$. Choose vectors
$\,\epsilon_i = (\epsilon_{i1},\ldots, \epsilon_{in}) \in {\bf R}^n \,$
for $1 \leq i \leq r$  such that, for all $i$ and all $s \not= t$, the numbers
$a_{is} + \epsilon_{is}$ and $a_{it} + \epsilon_{it}$ are distinct,
and  $$a_{is} + \epsilon_{is} < a_{it} + \epsilon_{it}
\quad \hbox{implies } \quad a_{is} \leq a_{it}.$$
The last condition is satisfied for all sufficiently small 
positive $\epsilon_i$.
Each vector $\epsilon_i$
defines a monomial $\,{\bf x}^{\epsilon_i}
\,= \, x_1^{\epsilon_{i1}} \cdots x_n^{\epsilon_{in}}\,$
with real exponents. We formally introduce the generic
monomial ideal (in a polynomial ring with real exponents):
$$\, M_\epsilon \quad := \quad \langle
\, m_1 \cdot {\bf x}^{\epsilon_1},\,
m_2 \cdot {\bf x}^{\epsilon_2},\, \ldots ,
m_r \cdot {\bf x}^{\epsilon_r}\, \rangle . $$
We call $M_{\epsilon}$ a {\it generic deformation} of $M$.
We abbreviate $\,\epsilon := (\epsilon_1, \ldots,\epsilon_r)$.
Let  ${\Delta}_{M_{\epsilon}}$ be the Scarf complex of $M_{\epsilon}$.
We now label the vertex of ${\Delta}_{M_{\epsilon}}$
corresponding to $m_i \cdot {\bf x}^{\epsilon_i}$
with the original monomial $m_i$.
Let ${\bf F}_\epsilon$ be the complex of
$S$-modules defined by this labeling of $\Delta_{M_\epsilon}$
as in Construction 1.1. \Box

\proclaim Theorem~5.2.
The complex ${\bf F}_{\epsilon}$ is a free
resolution of $S/M$ over $S$.

\noindent {\sl Proof: } Fix a monomial $m$.
Let $J$ be the largest subset of  $\{1,\dots ,r\}$
such that $m_J$ divides $m$.
The following conditions are equivalent
for a subset $I$ of $\{1,\ldots,r\}$:
$$ m_I \,\,\hbox{divides} \,\, m  \,\,\, \iff \,\,\,
 I \subseteq J \,\,\, \iff \,\,\,
 m_I \,\,\hbox{divides} \,\, m_J  \,\,\, \iff \,\,\,
 m_{I}(\epsilon) \,\,\hbox{divides} \,\, m_{J}(\epsilon). $$
Here $\,m_{I}(\epsilon) := lcm \bigl(
m_i {\bf x}^{\epsilon_i} : i \in I \bigr) $.
The last equivalence follows from our choice of the $\epsilon_{ij}$.
The  set of all faces of $\Delta_{M_\epsilon}$ which satisfy 
the four equivalent conditions above
is a contractible simplicial complex, by Theorem 4.1
applied to $M_\epsilon [m_{J}(\epsilon)]$.
Now apply Lemma 2.1 to $M$ and $m$
with $\Delta = \Delta_{M_\epsilon}$.  \Box

\proclaim Corollary~5.3. The Betti numbers of $M$ are
less or equal to those of any degeneration $M_\epsilon$,
that is, less or equal to the face numbers of the Scarf complex
$\Delta_{M_\epsilon}$.

\vskip .2cm

We emphasize that the Betti numbers of $M_{\epsilon}$ depend
on the choice of the generic deformation.
There are finitely many complexes $\Delta_{M_\epsilon}$
which can be obtained by Construction 5.1.
Each of them corresponds to $\epsilon$ lying in an
open cone in ${\bf R}^{r\cdot n}$.

\vskip .2cm

\noindent {\bf Example~5.4.} Consider Backelin's ideal
$M=\langle x^2,xy^2z,y^2z^2,yz^2w,w^2\rangle$
from Example~3.5. A generic deformation is
$M_\epsilon=\langle x^2,xy^2z,y^3z^3,yz^2w,w^2\rangle$.
Label the generators as $1,2,3,4,5$ in the given order.
The Scarf complex  of $M_\epsilon$
consists of the tetrahedron $\{1,2,4,5\}$ and the triangle $\{2,3,4\}$.
Relabeling its vertices by the generators of $M$ and applying Construction~1.1
we obtain a non-minimal free resolution ${\bf F}_{\epsilon}$ of $S/M$.
The Betti numbers of $S/M$ are $1,5,7,4,1$ while
the Betti numbers of $S/M_\epsilon$ are $1,5,8,5,1$.
Thus  ${\bf F}_\epsilon$ differs from the minimal resolution by a single
summand  $\, 0\, \rightarrow \, S\, \rightarrow \, S\, \rightarrow \, 0,$
placed in homological degrees $2$ and $3$.
However, this makes a big difference structurally: by Corollary~3.4,
the resolution
${\bf F}_{\epsilon}$ is a DG-algebra (with a simple multiplication rule)
while the minimal free resolution admits no DG-algebra structure at all.
Note that Taylor's resolution is one step longer than
${\bf F}_\epsilon$. It has Betti numbers $1,5,10,10,5,1$. \Box

\noindent {\bf Remark~5.5.}
While
choosing a degeneration of a
 monomial ideal we are guided
by the intuition of real exponent deformations.
In fact  no algebra in
a polynomial ring with real exponents is used. One can treat the
 $\epsilon_{ij}$ as symbolic quantities whose sole purpose
is to break ties whenever a variable has the
same non-zero degree in two distinct generators of $M$.
The Scarf complex of a generic monomial ideal depends
only on the order of the  generating exponents coordinatewise. Given an
 ideal with real generating exponents,  we can relabel the exponents
in each variable as integers while preserving their order. We will obtain a
monomial ideal  with integer exponents and the same Scarf complex. \Box

\noindent {\bf Example~5.6.}
A simple method for degenerating $M$ is to
slightly increase the exponents of the generators
$m_i$ in the order they are listed. For instance, we
can choose a positive integer $\nu > r$, and then set
$$ M_\epsilon \quad := \quad
\langle \, m_i \cdot (x_1 x_2 \cdots x_n)^{i/\nu}\,\,: \,\,
i = 1,2\ldots, r \, \rangle .$$
One can  avoid working with fractional exponents by taking instead
$$ M_\epsilon^{(\nu)} \quad := \quad
\langle \, m_i^\nu \cdot (x_1 x_2 \cdots x_n)^{i}\,\,: \,\,
i = 1,2\ldots, r \, \rangle .$$
Clearly, the generic monomial ideals
$ M_\epsilon $ and $M_\epsilon^{(\nu)}$
have the  same Scarf complex.

If $M$ is squarefree then
$\langle m_1,m_2^2,\dots ,m_r^{r} \rangle $
is a generic deformation of $M$. \Box

\beginsection  6. Bounds on  Betti numbers

In this section we address the following  problem:

\proclaim Upper Bound Problem.
Determine the maximal $i$-th Betti number $\beta_i(n,r) :=
\beta_i(M)$ among all
monomial ideals $M$ with $r$ minimal generators in $k[x_1,\ldots,x_n]$.

A related result in [Bi] and [Hu]
states that among all monomial ideals with fixed Hilbert function
the lexicographic ideal has maximal Betti numbers. In our problem
we do not fix the Hilbert function but only the number of generators.
The ideals attaining $\beta_i(n,r)$
are generally far from lexicographic; see (6.1) for an example.
It follows from Corollary 5.3 that  $\beta_i(n,r)$ is
attained by a monomial ideal $M$ which is generic. 
We may assume that $M$ is artinian, by the
following lemma:

\proclaim Lemma 6.1.
Let $\,M \, = \, \langle m_1,m_2,\ldots,m_r \rangle \,$ be a
generic monomial ideal where $\,m_1 \, =\,
x_1^{i_1} x_2^{i_2} \cdots x_n^{i_n} \,$
and $i_1 > deg_{x_1}(m_j)$ for $j \geq 2$.
If $\,M' \, = \, \langle x_1^{i_1},m_2,\ldots,m_r \rangle \,$
then the Scarf complex $\,\Delta_M \,$ is a subcomplex of
the Scarf complex $\,\Delta_{M'} $.

\noindent {\sl Proof: } Straightforward from the definitions. \Box

\proclaim Corollary 6.2.
$\beta_i(n,r)$ equals the maximal number of $i$-faces of any Scarf complex
$\Delta_M$, where $M$ runs over all artinian
generic monomial ideals $M$ with $r$ generators
in $n$ variables.

\noindent {\sl Proof: }
Apply Lemma 6.1 repeatedly until all variables appear to some power.
Take the resulting generic artinian monomial ideal $M$
and apply Corollary 3.3 (1).
\Box

Each Scarf complex $\Delta_M$ considered in Corollary 6.2 is
the boundary of a simplicial $n$-polytope with 
at least one facet removed, by Corollary 4.5.
In fact, in order to attain $\beta_i(n,r)$ it is enough
to consider simplicial $n$-polytopes with exactly one facet removed.
The {\it Upper Bound Theorem for Convex Polytopes}
(cf.~[Zi, Thm.~8.23]) implies the following  result:

\proclaim  Theorem~6.3. The Betti numbers of monomial
ideals satisfy the inequalities of the Upper Bound Theorem for
Convex Polytopes.
More precisely, if
 $c_{i}(n,r)$ denotes the number of $i$-dimensional faces of
the  cyclic $n$-polytope with $r$ vertices, then
$$  \eqalign{ \beta_{i}(n,r) \quad & \leq \quad c_{i}(n,r)
\qquad \qquad \hbox{for $1 \leq i \leq n-2,$} \cr
\hbox{and} \,\,\, \quad  \beta_{n-1}(n,r) \quad & \leq
\quad c_{n-1}(n,r) -1 . }  $$

An explicit formula for $ c_{i}(n,r)$
is given in [Zi, \S 8]. For instance,
$ c_1(3,r) =  3r-6  , \,
   c_2(3,r)  = 2r-4\,\,$ and $ \,
 c_1(4,r) = { r \choose 2 } ,\,
c_2(4,r)  = r (r - 3) , \,
c_3(4,r)  =  r (r-3) /2  $.

The number $\beta_{n-1}(n,r)$ coincides with
the maximal number of socle elements
modulo any artinian ideal
generated by $r$ monomials in $n$ variables.
This number  was studied recently  in [Ag],
where the socle elements are called ``outside corners''.

\proclaim Theorem 6.4. {\rm (Agnarsson 1996) }
The inequalities in Theorem 6.3 are equalities
for $n \leq 3$ and
for $n=4, r \leq 12 $, but
they are strict inequalities for $n=4, r \geq 13 $.

\noindent {\sl Proof: }
If $n \leq 3$ then all simplicial $n$-polytopes
with $r$ vertices have the same $f$-vector.
For $n = 4$ and $r=12$ Agnarsson constructs the
following generic artinian monomial ideal:
$$
\eqalign{
\langle \, & a^9 ,\, b^9 , \, c^9 ,\,  d^9 ,\,
a^6 b^7 c^4 d, \,
a^5 b^8 c^3 d^2 ,\,
a^8 b^5 c^2 d^3 ,\,
a^7 b^6 c d^4 , \cr & 
a^2 b^3 c^8 d^5 , \,\,
a b^4 c^7 d^6 ,\,\,
a^4 b c^6 d^7 ,\,\,
a^3 b^2 c^5 d^8 \, \rangle  . \cr} \eqno (6.1) $$
The Scarf complex of this ideal is a triangulation of the
tetrahedron with eight interior vertices.
It is {\it neighborly}, i.e., any two vertices are
connected by an edge.
Hence (6.1) has $\, c_1(4,12) = 66 \, $ first syzygies,
$\,c_2(4,12)  = 108 \,$ second syzygies and
$\,c_3(4,12)-1  = 53  \,$ third syzygies.
Taking subsets of the twelve generators, we obtain the assertion for
$n=4$ and $r \leq 12$.  In [Ag] it is proved that $\,\beta_1 (4,13) = 77$.
But $\,c_1(4,13) = 78$.
This implies $ \,\beta_2(4,13) <  c_2(4,13) \,$ and
$\,\beta_3(4,13) <  c_3(4,13) -1 \,$ by the Euler and
Dehn-Sommerville equations. \Box

Next we translate our problem into the
language of partially ordered sets.
We recall (for instance, from [Re] or [BT])
that the {\it order dimension} $odim({\cal P})$
of a finite poset ${\cal P}$ is the
smallest number $s$ of linear extensions $L_1,\ldots,L_s$ of ${\cal P}$
such that $\,L_1 \cap \dots \cap L_s \,= \,{\cal P} $.
If $\Delta$ is a  simplicial complex then  $odim(\Delta)$
denotes the order dimension of the face poset of $\Delta$.
It is well-known (see p.~59 in [Re])
that $\, odim(\Delta) \geq dim(\Delta) + 1$.
The case of equality is of special interest for us:

\proclaim Theorem 6.5.
\item{(a)} A simplicial complex $\Delta$ satisfies
$\,odim(\Delta) \leq n \,$ if and only if $\Delta$
is a subcomplex of the Scarf complex $\Delta_M$
for some generic monomial ideal $M$ in $k[x_1,\ldots,x_n]$.
\item{(b)} Let $\Delta $ be a triangulation of 
the $n$-ball whose boundary equals the boundary of an $(n-1)$-simplex.
Then $\,odim(\Delta)= n \, $ if and only if
$\Delta$ equals the Scarf complex $\Delta_M$ of a generic artinian monomial
ideal $M$ in $k[x_1,\ldots,x_n]$.

\noindent {\sl Proof: }
The if-direction in both (a) and (b) is seen as follows:
For $i=1,\ldots,n \,$ let $\, L_i$ denote the
linear extension of the face poset of $\Delta_M$ defined by
$\,\, \Delta_M \rightarrow {\bf R}, \,\, I \,\mapsto \,deg_{x_i}(m_I)
+ \epsilon {|I|} \,$, where $\epsilon $ is a small positive real.
The face poset of $\Delta_M$ coincides with   $\,L_1 \cap \dots \cap L_n$.

We next prove the only-if direction in (a).
Let $\Delta$ be a simplicial complex on $\{1,\ldots,r\}$
of order dimension at most $n$. Fix an embedding of
posets $\,\phi : \Delta \rightarrow {\bf N}^n \,$ such that
each coordinate of $\phi$ is a linear extension of $\Delta$.
We define the monomial ideal
$\, M \, = \, \langle m_1, m_2, \ldots,m_r \rangle \,$
where $\,m_i \, := \, {\bf x}^{\phi( \{ i \} )} $.
Let $I$ be any face of $\Delta$.  Note that
$\,m_I = lcm( m_i : i\in I)\,$ divides ${\bf x}^{\phi(I)}$.
We must show that
$I$ is a face of $\Delta_M$. Suppose not. Then there
exists a subset $J $ of $\{1,\ldots,r\}$
with $I \not= J$ but $m_I = m_J $.
If $J$ is not a subset of $I$ then
pick any $\,j \in J \backslash I $.
Then  $m_j = {\bf x}^{\phi({j})}$ does not divide
${\bf x}^{\phi(I)}$ (since $\phi$ is a poset embedding),
but $m_j$ does divide $m_I = m_J$, a contradiction.
If $J$ is a subset of $I$ then
$J$ is a proper face of $I$ in $\Delta$.
In this case we pick any $\,i \in I \backslash J $.
Now  $m_i$ does not divide
${\bf x}^{\phi(J)}$ (since $\phi$ is a poset embedding),
but $m_i$ does divide $m_I = m_J$, a contradiction.

We finally prove only-if in (b).
By part (a) there is a generic monomial ideal $M$
such that $\Delta$ is a subcomplex of $\Delta_{M}$.
We may assume that $\Delta $ and $\Delta_{M}$
have the same vertices. Applying Lemma 6.1
to the vertex labels of the $(n-1)$-simplex triangulated
by $\Delta$, we may also assume that $M$ is artinian.
Both $\Delta$ and $\Delta_M$ are triangulations
of the same $(n-1)$-simplex and $\Delta \subseteq \Delta_M$.
Hence $\Delta = \Delta_M$ . \Box

Theorem 6.5 implies that our
original problem can be rephrased as follows:

\proclaim Upper Bound Problem.
Determine the maximal cardinality $\beta_{i-1}(n,r)$ of a family ${\cal F}$
of $i$-sets in $\{1,\ldots,r\}$ such that
the simplicial complex spanned by ${\cal F}$
has order dimension $n$.

It was shown in 
[Sp, Theorem 2] (see also [Re, p.~61]) that
the complete graph on $r$ vertices has order dimension at least
$\, 2 + log_2 \bigl( log_2 ( r-1 ) \bigr) $.
Therefore  $\, \beta_1 (n,r) \, \leq \, { r \choose 2 } - 1 \,$
whenever $\,r \gg n $. Thus the inequalities in Theorem 6.2 are
all strict for $\, r \gg n \geq 4 $.
This extends Proposition 6.4.
For the reader's convenience
we give an algebraic reformulation and a proof of Spencer's result.

\proclaim Proposition 6.6. {\rm (Spencer 1971)  }
\hfill \break
Let $ M \subset k[x_1,\ldots,x_n]$ be an ideal
generated by $\, r > 2^{2^{n-1}} + 1$ monomials.
Then two of the minimal generators of $M$ are not linked by
a minimal first syzygy.

\noindent {\sl Proof: }
We use the classical theorem of Erd\"os and Szekeres that
if $m^2+1$ elements are ordered in two different ways
then some $(m+1)$-set is monotone under both orders.
By a simple induction we find that
if $\, 2^{2^{n-1}} + 1 \,$ elements are ordered in $n$
ways  then some triple $(m_1,m_2,m_3)$ is monotone under
all $n$ orders. If the elements are the minimal generators of a
monomial ideal in $k[x_1, \ldots,x_n]$ and the $n$ orders
are the degrees in $x_1,\ldots,x_n$ then this means that $m_2$ divides
$\, lcm(m_1,m_3)$, and we conclude that the pair
$(m_1,m_3)$ is not linked by a  minimal first syzygy
(i.e., $\{ m_1,m_3 \}$ is not an edge of the Scarf complex). \Box

\beginsection 7. Realizability

A central question in discrete geometry
is the {\sl Steinitz Problem} of characterizing
the face lattices of  all convex polytopes.
While this problem has a beautiful solution
due to Steinitz in dimension $3\,$ [Zi, Theorem 4.1], the
recent work of Richter-Gebert (see [RZ]) shows that
the Steinitz Problem for $4$-dimensional polytopes
is essentially equivalent to classifying all semi-algebraic
varieties (and hence intractible).
In this section we introduce a variant
of Steinitz' problem for monomial ideals.

\proclaim Realizability Problem.
Classify all simplicial complexes which are Scarf
complexes of generic (artinian) monomial ideals
with $r$ minimal generators in $k[x_1,\ldots,x_n]$.

Theorem 6.3 implies that this
can be rephrased in terms of order dimension:

\proclaim Realizability Problem.
Classify all (regular) triangulations of the $(n-1)$-simplex
whose face poset has order dimension $n$.

This problem is trivial for $n \leq 2$ variables.
In three variables it is non-trivial but completely
solved by results in [Sch] and [BT].
Theorem 6.1 in [BT] asserts that
every triangulation of the triangle has order dimension $3$.

\proclaim Theorem 7.1. {\rm (Schnyder 1989)}
For any triangulation  $\Delta$ of a triangle
there exists a generic artinian monomial ideal
$M \subset k[x,y,z]$ such that $\Delta = \Delta_M$.

The analogous result does not hold for $n=4$ and $r \geq 7 $.

\proclaim Theorem 7.2.
There exists a triangulation of the tetrahedron
with seven vertices which is not the Scarf complex
of any monomial ideal in four variables.

For the proof
we need one lemma.
Let $\Delta$ be any simplicial complex on
$\{1,\ldots,n, \ldots,r \}$
which is a triangulation of the $(n \!- \!1)$-simplex
$\{1,\ldots,n\}$.
A {\it labeling of $\Delta$} is a family of bijections
$\,\bigl\{ \,\phi_I \,: \, I  \rightarrow
\{x_1,\ldots,x_n\}  \bigr\}_{I \,
{\rm facet} \,{\rm of} \,\Delta }\,$ which satisfies the following two axioms:
\item{(A)} If $I$ is a facet of $\Delta \,$ and
$\, i \in I \,\cap \,\{1,\ldots,n\}$ then
$\phi_I (i) = x_i$.
\item{(B)}
If $\,I\,$ and $\,J \,$ are facets of $\Delta$
which share a common ridge, i.e.~$ I \cap J$ has
cardinality $n-1$, then  $\,\{ \,j \in I \cap  J \,: \,
\phi_I (j) \, =\,\phi_{J'} (j)  \}\,$
has cardinality $n-2$.

\proclaim Lemma 7.3.
The Scarf complex $\Delta_M$ of any
generic artinian monomial ideal $M$ possesses a
labeling $\phi^M$. We call the labelings
of the form $\phi^M$ {\it realizable}.

\noindent {\sl Proof:}
The labeling $\phi^M = \{ \phi^M_I \}_{I \, {\rm facet}}$
of $\Delta_M$ is defined as follows:
If $I $ is a facet of $\Delta_M$ and $j \in I $
then $\,\phi^M_I (j) \,$ is the unique variable
missing in $\,m_I /m_j $. The map
$\,\phi^M_I  : I \rightarrow \{x_1,\ldots,x_n\} \,$ is
a bijection, and it is straightforward to check the axioms (A) and (B).
\Box

\noindent {\sl Proof of Theorem 7.2: } Consider the following
triangulation of a tetrahedron $1234$:
$$ \Delta \,\, = \,\, \bigl\{  1237, 1245, 1256, 1267, 1347,
1457, 1567, 2345, 2356, 2367, 3456, 3467, 4567  \bigr\}.  $$
The complex $\Delta$ is the boundary of the
cyclic polytope $C_4(7) $ minus the facet $1234$.

By Lemma 7.3 it suffices to prove
that $\Delta$ does not admit any labeling.
Suppose on the contrary that $\Delta$ possesses a labeling.
Writing $a,b,c,d$ for the four variables, the following labels are
forced by the axioms:
$$
\vcenter{\halign{
&\hfil$#$\hfil\cr
\left[\matrix{ 2 \, 3 \, 4 \, 5 \cr b \, c \, d \, a \cr}\right]
&
\longrightarrow
&
\left[\matrix{ 3 \, 4 \, 5 \, 6 \cr c \, d \, b \, a \cr}\right]
&
\longrightarrow
&
\left[\matrix{ 3 \, 4 \, 6 \, 7 \cr c \, d \, b \, a \cr}\right]
\cr
& & & \searrow & \downarrow \cr
& & & &
\left[\matrix{ 4 \, 5 \, 6 \, 7 \cr d \, b \, c \, a \cr}\right]
&
\longrightarrow
&
\left[\matrix{ 1 \, 5 \, 6 \, 7 \cr a \, b \, c \, d \cr}\right]
\cr
}}
$$
$2345$ is labeled by axiom (A), and $3456$ and $3467$ successively by axioms
(A) and the ``ridge axiom'' (B). These last two tetrahedra each allow two
labelings of $4567$; the one shown is the common allowable labeling,
forcing the labeling shown on $1567$.

On the other hand, the following labels are also forced:
$$
\vcenter{\halign{
& $#$ \cr
\left[\matrix{ 1 \, 2 \, 3 \, 7 \cr a \, b \, c \, d \cr}\right]
&
\longrightarrow
&
\left[\matrix{ 1 \, 2 \, 6 \, 7 \cr a \, b \, d \, c \cr}\right]
\cr
}}
$$
Now, there is a contradiction to axiom (B) in the labelings
of $1267$ and $1567$. \Box

A monomial ideal
$\, M = \langle m_1,\ldots,m_r \rangle \subset k[x_1,\ldots,x_n] \,$
may be regarded as a non-negative integer point in the matrix space
${\bf R}^{r \times n}$. Two generic artinian monomial ideals $M$ and $M'$ are
called {\it equivalent} if $\,\Delta_{M} =\Delta_{M'} $.
We call $M$ and $M'$ {\it strongly equivalent} if
$\,\Delta_{M} =\Delta_{M'} \,$ and $\,\phi^{M} =\phi^{M'} \,$
(as defined in Lemma 7.3).
The equivalence class of $M$ is the union of
finitely many strong equivalence classes, one for
each realizable labeling of $\Delta_M$:

\vskip .2 cm

\noindent {\bf Remark 7.4.}
Each strong equivalence class consists of the interior integer
points in a convex polyhedral cone in the matrix space ${\bf R}^{r \times n}$.
\Box

\noindent {\bf Example 7.5.}
The  ideals
$\, M \, = \,\langle
x^4,
y^4,
z^4,
x y^2 z^3,
x^3 y z^2,
x^2 y^3 z
\rangle $ and $ M' \, = \,\langle
x^4,\allowbreak
y^4,\allowbreak
z^4,\allowbreak
x y^3 z^2,
x^2 y z^3,
x^3 y^2 z
\rangle \,  $
are equivalent but not strongly equivalent. Labeling the generators
$1,2,3,4,5,6$, their common Scarf complex is the boundary of an
octahedron minus a facet  $123$. It is depicted in Example 4.7.
This complex admits two realizable labelings which differ
precisely on the antipodal
facet $456$; the labeling on the left corresponds to $M$, and the labeling
on the right corresponds to $M'$:

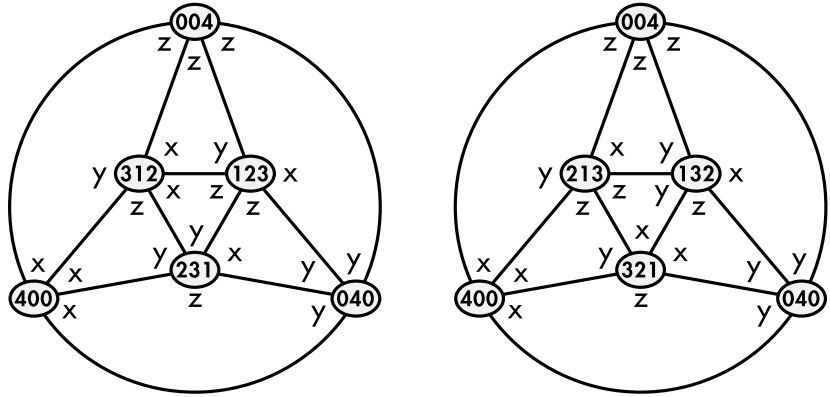
\draw4{309}{figD} 

Note how the exponents determine the labeling: Within each cell, each
variable labels the vertex having the largest exponent in that coordinate.
Conversely, a realizable labeling determines a strong equivalence class of
ideals: The labels of each variable can be thought of as the heads of
arrows forming a directed graph, giving the inequalities that exponents in
that coordinate must satisfy. The following directed graphs for $x$, $y$,
$z$ determine the class of ideals strongly equivalent to $M$:

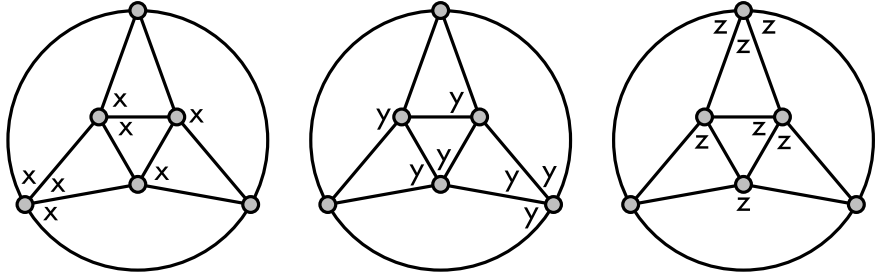
\draw1{328}{figE} 

\beginsection  8. Irreducible decomposition

We shall describe the irreducible
decomposition of a generic monomial ideal
$M$, that is, the unique minimal
expression of $M$ as an intersection of ideals of the form
$\, \langle x_{i_1}^{d_1}, x_{i_2}^{d_2}, \ldots, x_{i_j}^{d_j} \rangle $.
Choose an integer $D$ larger than the degree of any minimal
generator of $M$. We replace $M$ by the artinian ideal
$$ M^* \quad := \quad M \,\, + \,\, \langle
\, x_1^D, \, \, x_2^D, \, \ldots, \,
\, x_n^D\,\rangle. $$
Let $\Delta_{M^*}$ be the Scarf complex of $M^*$. This is
a pure $n$-dimensional simplicial complex on
$\{1,2,\ldots,r,r \!+ \! 1,\ldots,r \! + \! n\}$, where the index
$\,r \!+\!i \,$
is associated with the generator $x_i^D$.
If $I$ is a facet of $\Delta_{M^*}$ then we form the irreducible ideal
$$ M_I \quad := \quad \langle \,\, x_s^{p_s} \,\,:\,\,\,
 p_s = deg_{x_s}(m_I) \,\hbox{  and  }\, p_s < D \,\rangle. $$
Note that $M_I$ is independent of the choice of $D$ and
may have less than $n$ generators.

\proclaim Theorem 8.1. A generic monomial ideal $M$ is the
intersection of the irreducible ideals $M_I$,
where $I$ runs over all facets of the Scarf complex
$\Delta_{M^*}$. This intersection is irredundant.

\noindent {\sl Proof: }
We first show that $M$ is contained in $M_I$ for every
facet $I $ of $\Delta_{M^*}$. Let $m_j$ be any minimal
generator of $M$. If $j \in I $ then there exists a variable $x_s$
with $\, deg_{x_s}(m_I)  \, = \, deg_{x_s}(m_j) \, < \, D$.
This implies $m_j  \in M_I$.
If $j \not\in I $ then there is a variable $x_s$
with $ \, deg_{x_s}(m_I) \, < \, deg_{x_s}(m_j) \, < \, D $,
which implies $m_j  \in M_I$ as well.

For the reverse inclusion $\,\bigcap_I M_I \subseteq M$,
we shall prove that every $M$-standard monomial is
$M_I$-standard for some facet $I$ of $\Delta_{M^*}$.
Let $m $ be an $M$-standard monomial, i.e., $m \not\in M$.
We may choose $D \gg 0$ so that $m \not\in M^*$. Next we select
a monomial $\tilde m$ such that $\,m \cdot \tilde m \not\in M^*$
but $\,x_i \cdot m \cdot \tilde m \in M^* \,$ for $i = 1,2,\ldots,n$.
There exist unique (and necessarily distinct) minimal
generators $\,m_{j_1},\ldots,m_{j_n}$ of $M^*$ with the property
that $m_{j_i} $ divides $\,x_i \cdot m \cdot \tilde m $.
Setting $\, I \, := \, \{ j_1,\ldots,j_n \}$, we have
$\,x_1 x_2 \cdots x_n \cdot  m \cdot \tilde m  \, = m_I $.
This implies $\,m \not\in M_I $.

Finally, we must show that the intersection $\cap M_I$ over
all facets $I$ of $\Delta_{M^*}$ is irredundant.
Fix a facet $I$ and consider the monomial
$\,m \,:= \, m_I/(x_1 x_2 \cdots x_n)$.
Clearly, $m \not\in M_I$. It suffices to show that
$m \in M_J$ for all other facets $J$ of $\Delta_{M^*}$.
Fix another facet $J$. There exists an index $i \in I \backslash J $
and a variable $x_s$ such that
$ \,  deg_{x_s}(m_I) \geq  deg_{x_s}(m_i) \, > \,
deg_{x_s}(m_J) $. This implies
$ \, D \geq deg_{x_s}(m_I)
>  deg_{x_s}(m) \geq deg_{x_s}(m_J) \,$
and therefore $\,m \in M_J$. \Box

\proclaim Example 8.2. \rm
The seven irreducible components of the generic monomial ideal
$$ \langle \,
x y^2 z^3 , \, x^3 y z^2 , \, x^2 y^3 z  \, \rangle
\; = \;
\langle x \rangle \,\cap \,
\langle y \rangle \,\cap \,
\langle z \rangle \,\cap \,
\langle x^3, y^2 \rangle \,\cap \,
\langle y^3, z^2 \rangle \,\cap \,
\langle z^3, x^2 \rangle \,\cap \,
\langle x^3, y^3, z^3 \rangle  $$
 correspond to the seven triangles in the Scarf complex of
the artinian ideal
$ \,\langle x^4,\, y^4,\,\allowbreak z^4,\,\allowbreak x y^2
z^3,\,\allowbreak x^3 y z^2,\,\allowbreak x^2 y^3 z \,\rangle $. This Scarf
complex is depicted in Example 4.7. \Box

Suppose the monomial ideal $M$ is not generic.
Then we choose any degeneration  $M'$ as in Section~5
and use  Theorem 8.1 to compute the
irreducible decomposition of $M'$. This decomposition specializes
to a (generally not minimal) irreducible decomposition of $M$.

\proclaim Corollary 8.3.
A generic monomial ideal $M$ is Cohen-Macaulay if and only if it is pure.

\noindent {\sl Proof: }
The ideal $M$ being pure means that all irreducible components
$M_I$ have the same dimension. Each facet $I$ of $\Delta_{M^*}$
defines a component $M_I$ with
$$ dim(M_I) \quad = \quad  | \, I \cap \{r+1,\ldots,r+n\} \, | . $$
To compute the depth of $S/M$ we consider the
Scarf complex $\Delta_M$. Note that $\Delta_M$
is a subcomplex of $\Delta_{M^*}$
which may have dimension less than $n-1$ and is generally not pure.
The depth of $M$ equals the minimum of the numbers
$\,n - | J | \,$ where $J$ runs over all
facets of  $\Delta_{M}$. Every facet $J$ of
$\Delta_{M}$ extends to a facet $I$ of $\Delta_{M^*}$,
and, conversely, if $I$ is a facet of $\Delta_{M^*}$
then $\,I \cap \{1,\ldots,n\}$ is a face of $\Delta_{M}$. Therefore
$$ depth(S/M) \quad = \quad min \biggl\{
 | \, I \cap \{r+1,\ldots,r+n\} \, | \,\, : \,\,
\hbox{$I$ facet of $\Delta_{M^*}$} \,\,\biggr\}. $$
This proves that $M$ is pure if and only if $dim(M) = depth(M)$. \Box

\vskip .5cm

\noindent{\bf Acknowledgments.} We are  grateful to Herbert Scarf for
inspiring discussions during his visit to Berkeley in April 1996.
They got us started on this project. Many thanks to G\"unter Ziegler for
a very helpful discussion on the material in Section 7.
Dave Bayer and  Bernd Sturmfels
are partially supported by the National Science Foundation.
Irena Peeva and Bernd Sturmfels are partially supported
by the David and Lucile Packard Foundation.

\vskip 1.2cm

\centerline{\bf References}

\vskip .1cm

\item{[Ag]} Geir Agnarsson: On the number of outside corners
of monomial ideals, preprint, UC Berkeley 1996.

\item{[Av]} Luchezar Avramov: Obstructions to the existence of multiplicative
structures on minimal free resolutions,
{\sl Americal Journal of Mathematics} {\bf 103} (1981) 1--31.

\item{[Bi]} Anna Maria  Bigatti:
Upper bounds for the Betti numbers of a given Hilbert function,
{\sl Communications in Algebra} {\bf 21} (1993) 2317-2334.

\item{[BS]} Dave Bayer and Michael Stillman:
Computation of Hilbert functions,
{\sl Journal of Symbolic Computation} {\bf 14} (1992) 31--50.

\item{[BLSWZ]} Anders Bj\"orner, Michel Las Vergnas,
Bernd Sturmfels, Neil White and G\"unter Ziegler:
{\sl Oriented Matroids}, Cambridge University Press, 1993.

\item{[Bj]} Anders Bj\"orner: Topological Methods,
in: {\sl Handbook of Combinatorics}
(eds.\ R.~Graham, M.~Gr\"otschel, L.~Lov\'asz), North-Holland, (1995)
1819-1872.

\item{[BW]} Anders Bj\"orner and Michelle Wachs:
Shellings of non-pure complexes and posets,
{\sl Transactions of the American Mathematical Society}, to appear.

\item{[BT]} Graham Brightwell and
William Trotter: The order dimension of
convex polytopes, {\sl SIAM Journal of Discrete Mathematics}
{\bf 6} (1993) 230--245.

\item{[EK]} Shalom Eliahou and Michel Kervaire:
 Minimal resolutions of some monomial ideals,
{\sl Journal of Algebra} {\bf 129} (1990) 1--25.

\item{[Ho]} Melvin Hochster:
Cohen-Macaulay rings, combinatorics and simplicial complexes,
in {\sl Ring Theory II}
(B.R.~McDonald and R.~Morris, eds.),
Lecture Notes in Pure and Appl.~Math., No.~26,
Dekker, New York, (1977) 171--223.

\item{[Hu]} Heather Hulett: Maximum Betti numbers of homogeneous
ideals with a given Hilbert function,
{\sl Communications in Algebra} {\bf 21} (1993) 2335-2350.

\item{[Lee]} Carl Lee:  Regular triangulations of convex polytopes, in
{\sl Applied Geometry and Discrete Mathematics - The Victor Klee
Festschrift}, (P.~Gritzmann and B.~Sturmfels, Eds.), American
Math.~Soc., DIMACS Series {\bf 4}, Providence, R.I. (1991) 443--456.

\item{[Ly]} Gennady Lyubeznik:
A new explicit finite free resolution of ideals generated by monomials
in an $R$-sequence,
{\sl J.~Pure and Applied Algebra}
{\bf 51} (1988) 193-195.

\item{[Re]} Klaus Reuter:
On the order dimension of convex polytopes,
{\sl European Journal of Combinatorics}
{\bf 11} (1990) 57-63.

\item{[RZ]} J\"urgen Richter-Gebert
and G\"unter Ziegler: Realization spaces of $4$-polytopes are universal,
{\sl Bulletin of the American 
Math.~Society} {\bf 32} (1995) 403--412.

\item{[Sca]} Herbert Scarf: {\sl The Computation of Economic Equilibria},
Cowles Foundation Monograph {\bf 24}, Yale University Press, 1973.

\item{[Sch]}  Walter Schnyder:
Planar graphs and poset dimension,
{\sl Order} {\bf 5} (1989) 323--343.

\item{[Sp]} Joel Spencer:
Minimal scrambling sets of simple orders,
{\sl Acta Math.~Acad.~Scient.~Hungaricae}
{\bf 22} (1971) 349--353.

\item{[Sr]} Hema Srinivasan:
Algebra structures on some canonical resolutions,
{\sl Journal of Algebra}
{\bf 122} (1989) 150--187.

\item{[St]} Richard Stanley:
{\sl Combinatorics and Commutative Algebra}, Birkh\"auser, Boston, 1996.

\item{[Zi]} G\"unter Ziegler:
{\sl Lectures on Polytopes}, Springer Verlag, New York, 1994.

\vskip 1.5cm

Dave Bayer, Department of Mathematics,
Barnard College, Columbia University,
New York, NY 10027, USA,
{\tt bayer@math.columbia.edu}.

\vskip .5cm

Irena Peeva and Bernd Sturmfels,
Department of Mathematics,
University of California,
Berkeley, CA 94720, USA,
{\tt irena/bernd@math.berkeley.edu}.

\bye